\documentclass[12pt]{article}
\usepackage{amsbsy,amsfonts,amssymb}
\usepackage{rotating}
\usepackage{graphicx}

\newtheorem{theo}{Theorem}[section]
\newtheorem{lemma}[theo]{Lemma}
\newtheorem{prop}[theo]{Proposition}
\newtheorem{example}[theo]{Example}
\newtheorem{defi}[theo]{Definition}

\newtheorem{remark}[theo]{Remark}

\newcommand{\be}{\begin{equation}}
\newcommand{\ee}{\end{equation}}
\newcommand{\ba}{\begin{array}}
\newcommand{\ea}{\end{array}}
\newcommand{\dfrac}{\displaystyle \frac}
\newcommand{\m}{\frak m}

\begin{document}

\title{Genera of zero-divisor graphs with respect to ideals}

\author{
 Hsin-Ju Wang \thanks{e-mail:  hjwang@math.ccu.edu.tw} \\
Department of Mathematics, National Chung Cheng University,\\
          Chiayi 621, Taiwan}

\date{}

\maketitle

\noindent{\bf Abstract.} In this paper, we study the genera of
zero-divisor graphs with respect to ideals in finite rings.

\section{Introduction}
In this paper, $R$ will denote a commutative ring with $1$. Often,
but not always, we will also assume that $R$ is finite. If $S$ is
a subset of $R$, we denote $S-\{0\}$ by $S^*$. Also, we use
$\mathbb{N}$ for the natural numbers and $\mathbb{F}_q$ for the
finite field of $q$ elements. \par A subject of study linking
commutative ring theory with graph theory has been the concept of
the {\it zero-divisor graph} of a commutative ring. Let $R$ be a
ring. The zero-divisor graph of $R$, denoted $\Gamma(R)$, is a
simple graph whose vertices are the nonzero zero-divisors of $R$
with two distinct vertices $x$ and $y$ joined by an edge if and
only if $xy=0$. This definition was introduced by Anderson and
Livingston in \cite{al} and later was studied extensively in
\cite{amy}, \cite{am}, \cite{afll}, \cite{al}, \cite{s} and
\cite{wa}. Recently, Redmond extends the concept of zero-divisor
graphs to {\it zero-divisor graphs with respect to ideals}.
Namely, a simple graph $\Gamma_I(R)$ with vertices $\{x\in
R-I~|~xy\in I ~for~some~y\in R-I \}$ in which distinct vertices
$x$ and $y$ are adjacent if and only if $xy\in I$.

\par One of the primitive subjects of topological graph theory is to
embed a graph into a surface. In plain words, it is to draw a
graph on a surface so that there is no crossing for any two edges.
One simple question which one may ask is that "What kind of rings
can have $\gamma(\Gamma(R))=g$ (see definition in section 2)?".
For example, in \cite{cw}, Chiang-Hsieh and Wang find all finite
rings that have genus one. When concerning to  zero-divisor graphs
with respect to ideals, similar question also has been discussed.
In \cite{r}, Redmond find a sufficient and necessary condition for
which $\Gamma_I(R)$ is planar. The goal of this note is to find
all finite rings and ideals for which the genus of $\Gamma_I(R)$
is one.

\par In order to obtain our main results, we review in section 2 some background
from graph theory and derive a criterion for a graph to have
bigger genus than its subgraph. \par In section 3, we establish
the relationship between the graph $\Gamma(R/I)^{(t)}$ and the
graph $\Gamma_I(R)$ if $||I|=t$, where $G^{(t)}$ is defined in
section 2 for any finite simple graph. \par In section 4, we first
show that if $\omega(\Gamma(R/I))\geq 4$ or
$\gamma(\Gamma(R/I))\geq 1$, then $\gamma(\Gamma_I(R))\geq 2$.
Then we state and prove some equivalent conditions for which
$\gamma(\Gamma_I(R))=1$.

\section{Preliminary}
We review some background from graph theory in this section.
\par A simple graph $G$ is an ordered pair of disjoint sets $(V, E)$ such
that $V=V(G)$ is the vertex set of $G$ and $E=E(G)$ is the edge
set of $G$. The {\it order} of a graph $G$, written by $|G|$, is
the cardinality of $V(G)$. A {\it subgraph} of $G$ is a graph
having all of its vertices and edges in $G$. Let $V'\subseteq
V(G)$; then $G-V'$ is the subgraph of $G$ obtained by deleting the
vertices in $V'$ and all edges incident with them. Similarly, if
$E'\subseteq E(G)$, then $G-E'$ is the subgraph of $G$ obtained by
deleting the edges in $E'$. For any set $S$ of vertices of $G$, we
use the symbol $<S>$ for the subgraph induced by $S$.\par A {\it
bipartite} graph (bigraph) $G$ is a graph whose point set $V$ can
be partitioned into two subsets $V_1$ and $V_2$ such that every
line of $G$ joins $V_1$ with $V_2$. If $G$ contains every line
joining $V_1$ and $V_2$, then $G$ is a {\it complete bipartite}
graph. If $|V_1|=m$ and $|V_2|=n$, we use the symbol $K_{m,n}$ for
the complete bipartite graph. Moreover, if $m=1$ or $n=1$, then
$G$ is called a {\it star graph}. A graph in which each pair of
distinct vertices is joined by an edge is called a {\it complete
graph}. We use $K_n$ for the complete graph with $n$ vertices. For
a graph $G$, the {\it clique number}, $\omega(G)$, is the greatest
integer $n\geq 1$ such that $K_n$ is a subgraph of $G$. If $K_n$
is a subgraph of $G$ for every $n$ then $\omega(G)=\infty$.
\begin{remark}
\emph{In \cite[Theorem~7.2]{r}, the author states that
$\omega(G)\leq 2$ if and only if $G$ contains no cycles. However,
this statement is wrong as if $G=C_4$ (the cycle of length $4$)
satisfies $\omega(G)\leq 2$.}
\end{remark}
\par  A simple graph is said to be {\it planar} if it
can be drawn in the plane or on the surface of a sphere. It is
known that $K_{3,3}$ and $K_5$ are not planar and can be drawn
without crossings on the surface of a torus. The torus can be
thought of as a sphere with one {\it handle}. More generally, a
surface is said to be of genus $g$  if it is topologically
homeomorphic to a sphere with $g$ handles. Thus the genus of a
sphere is 0 and the one of torus is one. A graph can be drawn
without crossings on the surface of genus $g$, but not on one of
genus $g-1$, is called a graph of genus $g$. We write $\gamma(G)$
for the genus of a graph $G$. Therefore
$\gamma(K_{3,3})=\gamma(K_5)=1$. A well-known fact is that if $G$
is a connected graph of genus $g$, then any presentation of $G$ on
a surface of genus $g$ satisfies $n-m+f=2-2g$ and $2m\geq 3f$,
where $n, m, f$ are the vertices, edges and faces according to the
presentation. There are some known results from \cite{ry} which
will be used later.

\begin{lemma}
\label{genuskn} $\gamma(K_n)=\{\dfrac{1}{12}(n-3)(n-4)\}$, where
$\{x\}$ is the least integer that is greater than or equal to $x$.
In particular, $\gamma(K_n)=1$ if $n=5, 6, 7$.
\end{lemma}
\begin{proof}
See, for example \cite{ry}.
\end{proof}
\begin{lemma}
\label{genusbi} $\gamma(K_{m,n})=\{\dfrac{1}{4}(m-2)(n-2)\}$,
where $\{x\}$ is the least integer that is greater than or equal
to $x$. In particular, $\gamma(K_{4,4})=\gamma(K_{3,n})=1$ if
$n=3, 4, 5, 6$.
\end{lemma}
\begin{proof}
See, for example \cite{ry}.
\end{proof}

For later use, we introduce the following notion:
\begin{defi}
\label{maindefi1} Let $G$ be a finite simple graph with vertex set
$\{w_1, \dots, w_m\}$ and edge set $E(G)$; then $G^{(t)}$ is the
finite simple graph with vertex set $\{w_{ij}~|~1\leq i\leq
t,~1\leq j\leq m\}$ and edge set $\{w_{ij}w_{kl}~|~j\neq
l,~w_jw_l\in E(G) \}$.
\end{defi}
\begin{remark}
\emph{Observe that $w_{ij}w_{kl}$ is not an edge of $G^{(t)}$ if
$j=l$.}
\end{remark}
\begin{example}
\label{gexam1} If $G=K_{m,n}$ then $G^{(t)}=K_{mt,nt}$. If $G=K_n$
then $G^{(t)}=K_{t, \dots, t}$.
\end{example}

There are some results concerning the genera of $G^{(t)}$.

\begin{prop}
\label{gprop}  Let $G$ be a finite simple graph. Then the
following hold:
\begin{description}
\item{(a)} If $K_2$ is a subgraph of $G$, then
$\gamma(G^{(5)})\geq 2$. \item{(b)} If $K_5$ is a subgraph of $G$,
then $\gamma(G^{(2)})\geq 2$.  \item{(c)} If $K_{1,3}$ is a
subgraph of $G$, then $\gamma(G^{(3)})\geq 2$. \item{(d)} If
$K_{2,3}$ is a subgraph of $G$, then $\gamma(G^{(2)})\geq 2$.
\end{description}
\end{prop}
\begin{proof} (a) We may assume that $G=K_2$. Then by
Example~\ref{gexam1}, $G^{(5)}=K_{5,5}$. Therefore,
$\gamma(G^{(5)})\geq 2$ by Lemma~\ref{genusbi}.\\  (b) We may
assume that $G=K_5$. Let $V(G)=\{w_1, \dots, w_5\}$ and
$E(G)=\{w_iw_j~|~1\leq i< j\leq 5\}$. Then
$V(G^{(2)})=\{w_{ij}~|~1\leq i\leq 2,~1\leq j\leq 5\}$ and
$E(G^{(2)})=\{w_{ij}w_{kl}~|~j\neq l\}$. Let $u_1=w_{11}$,
$u_2=w_{12}$, $u_3=w_{21}$, $u_4=w_{22}$, $v_1=w_{13}$,
$v_2=w_{14}$, $v_3=w_{15}$, $v_4=w_{23}$, $v_5=w_{24}$ and
$v_6=w_{25}$; then $u_iv_j$ is an edge of $G^{(2)}$, so that
$K_{4,6}$ is a subgraph of $G^{(2)}$, it follows that
$\gamma(G^{(2)})\geq 2$ by Lemma~\ref{genusbi}.\\ (c) We may
assume that $G=K_{1,3}$. Then by Example~\ref{gexam1},
$G^{(3)}=K_{3,9}$. Therefore, $\gamma(G^{(3)})\geq 2$ by
Lemma~\ref{genusbi}.\\ (d) We may assume that $G=K_{2,3}$. Then by
Example~\ref{gexam1}, $G^{(2)}=K_{4,6}$. Therefore,
$\gamma(G^{(2)})\geq 2$ by
Lemma~\ref{genusbi}.\\
\end{proof}
To end this section, we illustrate the following useful Lemma.
\begin{lemma}
\label{keylem} Let $G$ be a connected graph and $G_1$, $G_2$ be
subgraphs of $G$. Suppose the following hold:\\
(a) $V(G)=V(G_1)\cup V(G_2)$ and $V(G_1)\cap V(G_2)=\emptyset$.\\
(b) $V(G_1)=\{v_1, \dots, v_4\}$ and $G_1=K_4$.\\
(c) For every $i$, there is an edge in $G$ joins $v_i$ with $G_2$.\\
Then $\gamma(G)>\gamma(G_2)$.
\end{lemma}
\begin{proof}
Let $g=\gamma(G_2)$. Suppose on the contrary that
$\gamma(G)=\gamma(G_2)$. Assume that $\{F_1, \dots, F_m\}$ and
$\{F'_1, \dots, F'_{m'}\}$ are faces of $G$ and $G_2$ when drawing
$G$ and $G_2$ on a surface $S_g$ of genus $g$. Then $\{F'_1,
\dots, F'_{m'}\}$ can be obtained by deleting $v_1, \dots, v_4$
and all edges incident with $v_1, \dots, v_4$ from the
presentation of $G$. Further, $\{F_1, \dots, F_m\}$ can be
recovered by inserting $v_1, \dots, v_4$ and all edges incident
with $v_1, \dots, v_4$ into the presentation corresponding to
$G_2$. Let $F'_{t_i}$ denote the face of $G_2$ into which $v_i$ is
inserted during the recovering process from $G_2$ to $G$. We note
that $v_iv_j\in E(G)$ for every $i\neq j$; therefore all $v_i$
should be inserted into the same face, say $F=F'_1$, of $G_2$ to
avoid any crossings, i.e., $t_1=\cdots=t_4=1$. Let $e_i$ be an
edge join $v_i$ to $G_2$ for $i=1, \dots, 4$. After inserting
$G_1$ into $F$ we obtain Figure 1 . However, it is easy to see
from Figure 1 that we can not insert $e_1, \dots, e_4$ into $F$
without crossings, a contradiction. Therefore, we conclude that
$\gamma(G)>\gamma(G_2)$.
\end{proof}

\begin{center}
\includegraphics{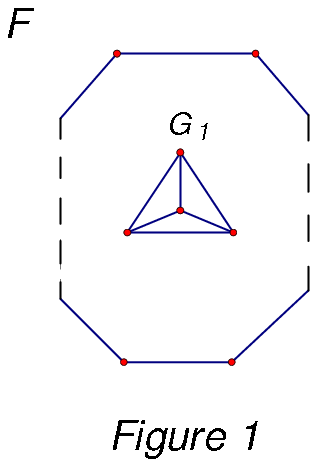}
\end{center}

\section{Basis properties of $\Gamma_I(R)$}
Throughout, let $R$ be a commutative ring with $1$ and $I$ be an
ideal of $R$. We use $\Gamma(R)$ to denote the zero-divisor graph
of $R$. In \cite{r}, Redmond extends the concept of zero-divisor
graphs as follows:
\begin{defi}
\label{maindefi2} Let $R$ be a commutative ring and let $I$ be an
ideal of $R$. We define a simple graph $\Gamma_I(R)$ with vertices
$\{x\in R-I~|~xy\in I ~for~some~y\in R-I \}$, where distinct
vertices $x$ and $y$ are adjacent if and only if $xy\in I$.
\end{defi}

\begin{center}
\includegraphics{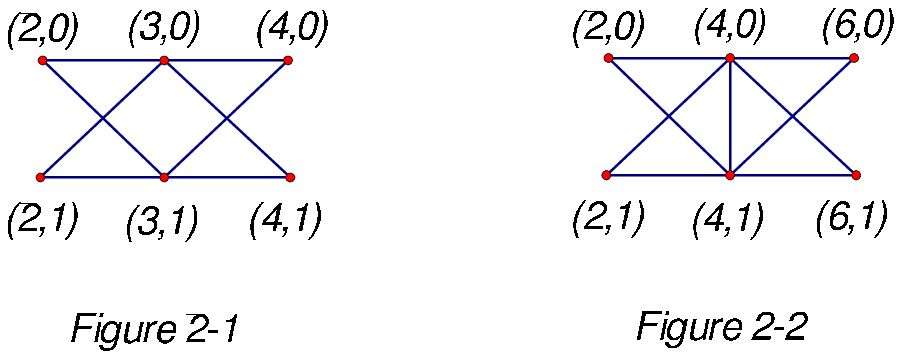}
\end{center}

\begin{example}
\emph{In Figure 2-1, $R=\mathbb{Z}_6\times \mathbb{Z}_2$ and
$I=0\times \mathbb{Z}_2$. In Figure 2-2, $R=\mathbb{Z}_8\times
\mathbb{Z}_2$ and $I=0\times \mathbb{Z}_2$.}
\end{example}

\begin{remark}
\emph{\label{keyrem1} (a) By Definition~\ref{maindefi2}, if
$\{x_{\lambda}+I~|~\lambda \in \Lambda\}$ is the set of nonzero
zero-divisors of $R/I$, then the vertex set of  $\Gamma_I(R)$ is
$\{x_{\lambda}+a~|~\lambda \in \Lambda,~a\in I\}$, therefore $|\Gamma_I(R)|=|I|\cdot |\Gamma(R/I)|$.  \\
(b) By Definition~\ref{maindefi1} and Definition~\ref{maindefi2},
if $|I|=t$ then $\Gamma(R/I)^{(t)}$ is a subgraph of
$\Gamma_I(R)$. More precisely,  if $\{x_1+I, \dots, x_m+I\}$ is
the set of nonzero zero-divisors of $R/I$, $I=\{a_1, \dots, a_t\}$
and $G=\Gamma(R/I)$, then the vertex set of $\Gamma_I(R)$ is
$\{x_j+a_i~|~1\leq i\leq t,~1\leq j\leq m\}$ and the edge set of
$\Gamma_I(R)$ is $\{(x_j+a_i)(x_l+a_k)~|~x_jx_l\in I\}$. However,
the edge set of $G^{(t)}$ is $\{(x_j+a_i)(x_l+a_k)~|~j\neq
l,~x_jx_l\in I\}$. Thus, $\Gamma(R/I)^{(t)}$ is a subgraph of
$\Gamma_I(R)$. Also, from this one can easily see that
$\Gamma(R/I)^{(t)}=\Gamma_I(R)$ if and only if $I$ is radical.\\
(c) From the above one can see that $\Gamma_I(R)$ is depend only
on $R/I$ and $|I|$. Therefore the graph $\Gamma_I(R)$ and the
graph $\Gamma_{0\times A}(R/I\times A)$ are the same, where $A$ is
any ring with $|A|=|I|$.}
\end{remark}

We list some properties of $\Gamma_I(R)$ in the following. Part of
them come from \cite{mpy} and \cite{r}.
\begin{theo}
$\Gamma_I(R)$ is finite if and only if $R$ is finite. Moreover,
$\Gamma_I(R)$ is connected.
\end{theo}

\begin{theo}
Let $R$ be a finite ring and $I$ be a proper nonzero ideal of $R$.
Then the diameter of $\Gamma_I(R)$ is less than or equal to $3$.
If $\Gamma_I(R)$ contains a cycle then $gr(\Gamma_I(R))\leq 4$.
\end{theo}

\begin{theo}
Let $I$ be a finite ideal of a ring $R$ and let $J$ be a finite
ideal of a ring $S$ such that $\sqrt{I}=I$ and $\sqrt{J}=J$. Then
the following hold: \begin{description} \item{(a)} If $|I|=|J|$
and $\Gamma(R/I)\cong \Gamma(S/J)$, then $\Gamma_I(R)\cong
\Gamma_J(S)$. \item{(b)} If $\Gamma_I(R)\cong \Gamma_J(S)$, then
$\Gamma(R/I)\cong \Gamma(S/J)$
\end{description}
\end{theo}

\begin{theo}
Let $I$ be a radical ideal of a ring $R$ such that $I=\cap_{1\leq
i\leq t} \frak{p}_i$ is a minimal primary decomposition, where
$\frak{p}$ is a prime ideal for every $i$. Then
$\omega(\Gamma_I(R))=n$.
\end{theo}

\section{The genera of $\Gamma_I(R)$}
In this section, we find some sufficient and necessary conditions
for a finite ring $R$ and a nonzero ideal $I$ of $R$ to have
$\gamma(\Gamma_I(R))\leq 1$. To begin with, we state some known
results and list some easy observations:
\begin{theo}
\label{rtheo} \cite[Theorem~7.2]{r} Let $R$ be a finite ring and
$I$ be a proper nonzero ideal of $R$ which is not a prime ideal.
Then $\Gamma_I(R)$ is planar if and only if $\Gamma(R/I)$ contains
no cycles and either (a) $|I|=2$ or (b) $|\Gamma(R/I)|=1$ and
$|I|\leq 4$.
\end{theo}

\begin{theo}
\label{watheo1} \cite{cw} Let $(R, \m)$ be a finite local ring
which is not a field. Then
$\Gamma(R)$ is planar if and only if $R$ is isomorphic to one of
the following 29 rings.

\medskip
\noindent
$\mathbb{Z}_4$, %
$\mathbb{Z}_8$, %
$\mathbb{Z}_9$, %
$\mathbb{Z}_{16}$, %
$\mathbb{Z}_{25}$, %
$\mathbb{Z}_{27}$, %
$\dfrac{\mathbb{Z}_2[x]}{(x^2)}$, %
$\dfrac{\mathbb{Z}_2[x]}{(x^3)}$, %
$\dfrac{\mathbb{Z}_2[x]}{(x^4)}$, %
$\dfrac{\mathbb{Z}_2[x,y]}{(x^2, xy, y^2)}$, %
$\dfrac{\mathbb{Z}_2[x,y]}{(x^2, y^2)}$, %
$\dfrac{\mathbb{Z}_2[x,y]}{(x^3, xy, y^2-x^2)}$, %
$\dfrac{\mathbb{F}_4[x]}{(x^2)}$, %
$\dfrac{\mathbb{Z}_3[x]}{(x^2)}$, %
$\dfrac{\mathbb{Z}_3[x]}{(x^3)}$, %
$\dfrac{\mathbb{Z}_4[x]}{(x^2)}$, %
$\dfrac{\mathbb{Z}_4[x]}{(x^2+x+1)}$, %
$\dfrac{\mathbb{Z}_4[x]}{(2x,x^2)}$, %
$\dfrac{\mathbb{Z}_4[x]}{(x^2-2,x^4)}$, %
$\dfrac{\mathbb{Z}_4[x]}{(x^3-2,x^4)}$, %
$\dfrac{\mathbb{Z}_4[x]}{(x^2-2,x^3)}$, %
$\dfrac{\mathbb{Z}_4[x]}{(x^3, x^2-2x)}$, %
$\dfrac{\mathbb{Z}_4[x]}{(x^3+x^2-2,x^4)}$, %
$\dfrac{\mathbb{Z}_4[x,y]}{(x^2, y^2, xy-2)}$, %
$\dfrac{\mathbb{Z}_4[x,y]}{(x^3,x^2-2, xy, y^2-2)}$, %
$\dfrac{\mathbb{Z}_5[x]}{(x^2)}$, %
$\dfrac{\mathbb{Z}_8[x]}{(x^2-4, 2x)}$, %
$\dfrac{\mathbb{Z}_9[x]}{(x^2-3,x^3)}$, %
$\dfrac{\mathbb{Z}_9[x]}{(x^2+3, x^3)}$. \\
\end{theo}

\begin{theo}
\label{watheo2} \cite{cw} Let $(R, \m)$ be a finite local ring
which is not a field. Then $\gamma(\Gamma(R))=1$ if and only if
$R$ is isomorphic to one of the following 18 rings.

\medskip \noindent $\mathbb{Z}_{32}$, %
$\mathbb{Z}_{49}$, %
$\dfrac{\mathbb{Z}_2[x]}{(x^5)}$, %
$\dfrac{\mathbb{F}_8[x]}{(x^2)}$, %
$\dfrac{\mathbb{Z}_2[x,y]}{(x^3, xy, y^2)}$, %
$\dfrac{\mathbb{Z}_2[x, y, z]}{(x, y, z)^2}$, %
$\dfrac{\mathbb{Z}_4[x]}{(x^3+x+1)}$, %
$\dfrac{\mathbb{Z}_4[x]}{(x^3-x+1)}$, %
$\dfrac{\mathbb{Z}_4[x]}{(x^3-2,x^5)}$, %
$\dfrac{\mathbb{Z}_4[x]}{(x^4-2, x^5)}$, %
$\dfrac{\mathbb{Z}_4[x]}{(x^4+x^3-2, x^5)}$,
$\dfrac{\mathbb{Z}_4[x]}{(x^3,2x)}$, %
$\dfrac{\mathbb{Z}_4[x,y]}{(x^3, x^2-2, xy, y^2)}$, %
$\dfrac{\mathbb{Z}_4[x,y]}{(2x ,2y, x^2,xy,y^2)}$, %
$\dfrac{\mathbb{Z}_7[x]}{(x^2)}$, %
$\dfrac{\mathbb{Z}_8[x]}{(x^2,2x)}$, %
$\dfrac{\mathbb{Z}_8[x]}{(x^2-2, x^5)}$, %
$\dfrac{\mathbb{Z}_8[x]}{(3x^2-2, x^5)}$. %
\end{theo}

\begin{theo} \label{watheo3} \cite{cw}  Let $R$ be a
finite ring which has exactly two maximal ideals; then $\Gamma(R)$
is planar if and only if $R$ is isomorphic to one of the following
15 rings.

\noindent
$\mathbb{Z}_2\times \mathbb{F}_q$, 
$\mathbb{Z}_3\times \mathbb{F}_q$, %
$\mathbb{Z}_2\times \mathbb{Z}_9$, %
$\mathbb{Z}_2\times \dfrac{\mathbb{Z}_3[x]}{(x^2)}$, %
$\mathbb{Z}_2\times \mathbb{Z}_4$, %
$\mathbb{Z}_2\times \dfrac{\mathbb{Z}_2[x]}{(x^2)}$, %
$\mathbb{Z}_2\times \dfrac{\mathbb{Z}_2[x]}{(x^3)}$, %
$\mathbb{Z}_2\times \dfrac{\mathbb{Z}_4[x]}{(x^2-2, x^3)}$, %
$\mathbb{Z}_2\times \mathbb{Z}_8$, %
$\mathbb{Z}_3\times \mathbb{Z}_9$, %
$\mathbb{Z}_3\times \dfrac{\mathbb{Z}_3[x]}{(x^2)}$, %
$\mathbb{Z}_3\times \mathbb{Z}_4$, %
$\mathbb{Z}_3\times \dfrac{\mathbb{Z}_2[x]}{(x^2)}$, %
$\mathbb{Z}_2\times \mathbb{Z}_2\times \mathbb{Z}_2$,
$\mathbb{Z}_2\times \mathbb{Z}_2\times \mathbb{Z}_3$.

\end{theo}

\begin{lemma}\label{mainlem1}
If $(R, \m)$ is finite local and $k$ is the smallest integer for
which $\m^k=0$, then $|\m^i|=t^{n_i}|\m^{i+1}|$ for $i=0, \cdots,
k-1$ if $t=|R/\m|$. In particular, $|R|=t^n$ for some $n$.
\end{lemma}
\begin{proof}
The assertion follows from the fact that $\m^i/\m^{i+1}$ is a
nonzero vector space over $R/\m$  for $i=0, \cdots, k-1$ and
$|R/\m|=t$.
\end{proof}
\begin{lemma}\label{mainlem2}
Let $(R, \m)$ be a finite local ring with $\m^2\neq 0$ and
$|R/\m|\geq  3$. Then $\gamma(\Gamma(R)^{(2)})\geq 2$.
\end{lemma}
\begin{proof}
Let $k$ be the smallest integer for which $\m^k=0$; then $k\geq 3$
by assumption. Let $t=|R/\m|$; then $t\geq 3$ by assumption.
Observe that $|\m^{k-1}-\{0\}|\geq t-1\geq 2$ and
$|\m^{k-2}-\m^{k-1}|\geq (t-1)|\m^{k-2}|\geq 6$ by
Lemma~\ref{mainlem1}. Therefore there are distinct elements $u_1,
u_2\in \m^{k-1}-\{0\}$ and distinct elements $v_1, \dots, v_6\in
\m^{k-2}-\m^{k-1}$. Since $k\geq 3$, $u_iv_j=0$ for every $i, j$.
Thus $K_{2,3}$ is a subgraph of $\Gamma(R)$ and then
$\gamma(\Gamma(R)^{(2)})\geq 2$ by Proposition~\ref{gprop}(d).
\end{proof}

\begin{lemma}\label{mainlem3}
Let $(R, \m)$ be a finite local ring such that $\m^2\neq 0$. If
$\Gamma(R)$ contains no cycles, then $|R/\m|=2$.
\end{lemma}
\begin{proof}
If $|R/\m|\geq 3$, then from the proof of Lemma~\ref{mainlem2}
$u_1-u_2-v_1-u_1$ is a triangle, a co
ntradiction. Thus,
$|R/\m|=2$.
\end{proof}

\begin{lemma}\label{mainlem4}
Let $R\cong \mathbb{Z}_2\times S$, where $S$ is a finite local
ring. Then the following hold: \begin{description} \item{(a)} If
$|\Gamma(S)|\leq 1$, then $\Gamma(R)$ is planar  and contains no
cycles. \item{(b)}  If $|\Gamma(S)|\geq 2$, then $K_3$ and
$K_{2,3}$ are subgraphs of $\Gamma(R)$.
\end{description}
\end{lemma}
\begin{proof}
(a) If $S$ is a finite field, then $\Gamma(R)$ is of course planar
and contains no cycles. By Theorem~\ref{watheo1}, if
$|\Gamma(S)|=1$, then $S\cong \mathbb{Z}_4$ or $S\cong
\dfrac{\mathbb{Z}_2[x]}{(x^2)}$, so that $R\cong
\mathbb{Z}_2\times \mathbb{Z}_4$ or $R\cong
 \mathbb{Z}_2\times \dfrac{\mathbb{Z}_2[x]}{(x^2)}$, it follows that $\Gamma(R)$ are trees by \cite[Figure~11]{r}. \\
(b) Suppose that $|\Gamma(S)|=2$. Then $S\cong \mathbb{Z}_9$ or
$S\cong \dfrac{\mathbb{Z}_3[x]}{(x^2)}$ by Theorem~\ref{watheo1}.
In either cases, there are two distinct nonzero zero-divisors $a,
b\in S$ such that $ab=a^2=b^2=0$. Let $u_1=(0, a)$, $u_2=(0, b)$,
$v_1=(1, 0)$, $v_2=(1, a)$ and $v_3=(1, b)$; then $\{u_1, u_2,
v_1, v_2, v_3\}$ are part of a vertex set of $\Gamma(R)$ such that
$u_iv_j=0$ for all $i, j$, so that $K_{2,3}$ is a subgraph of
$\Gamma(R)$. Moreover, $u_1-u_2-v_1-u_1$ is a triangle contained
in $\Gamma(R)$. \par Finally assume that $|\Gamma(S)|\geq 3$. Then
there are three distinct nonzero zero-divisors $a, b, c\in S$ such
that $ab=ac=0$. Let $u_1=(0, b)$, $u_2=(0, c)$, $v_1=(1, 0)$,
$v_2=(0, a)$ and $v_3=(1, a)$; then $\{u_1, u_2, v_1, v_2, v_3\}$
are part of a vertex set of $\Gamma(R)$ such that $u_iv_j=0$ for
all $i, j$, so that $K_{2,3}$ is a subgraph of $\Gamma(R)$.
Moreover, $v_1-v_2-u_1-v_1$ is a triangle contained in
$\Gamma(R)$.
\end{proof}

\begin{lemma}\label{mainlem5}
Let $(R, \m)$ be a finite local ring and $\Gamma(R)=K_3$. Then
$R\cong \dfrac{\mathbb{Z}_2[x,y]}{(x^2, xy, y^2)}$ or
$\dfrac{\mathbb{Z}_4[x]}{(2x, x^2)}$ or $R\cong
\dfrac{\mathbb{F}_4[x]}{(x^2)}$ or $R\cong
\dfrac{\mathbb{Z}_4[x]}{(x^2+x+1)}$. In particular, $\m^2=0$.
\end{lemma}
\begin{proof}
Among all rings in Theorem~\ref{watheo1}, there are exactly four
rings $\dfrac{\mathbb{Z}_2[x,y]}{(x^2, xy, y^2)}$,
$\dfrac{\mathbb{Z}_4[x]}{(2x, x^2)}$,
$\dfrac{\mathbb{F}_4[x]}{(x^2)}$ and
$\dfrac{\mathbb{Z}_4[x]}{(x^2+x+1)}$ satisfy the assumption.
Moreover, it is easy to see that $\m^2=0$ if $R$ is one of the
above three rings.
\end{proof}

\begin{lemma}\label{mainlem6}
Let $R\cong R_1\times \cdots \times R_k$, where $R_i$ is a finite
local ring for every $i$ and $k\geq 3$. If
$\gamma(\Gamma(R)^{(2)})\leq 1$, then $k=3$ and $R_i\cong
\mathbb{Z}_2$ for every $i$.
\end{lemma}
\begin{proof}
The first statement follows if we can show that $k=4$ implies that
$\gamma(\Gamma(R)^{(2)})\geq 2$. For this, suppose that $k= 4$.
Let $u_1=(1, 0, 0, 0)$, $u_2=(0, 1, 0, 0)$, $u_3=(1, 1, 0, 0)$,
$v_1=(0, 0, 1, 0)$, $v_2=(0, 0, 0, 1)$ and $v_3=(0, 0, 1, 1)$;
then $u_iv_j=0$ for all $i, j$, so that $K_{3,3}$ is a subgraph of
$\Gamma(R)$ it follows that $K_{6,6}$ is a subgraph of
$\Gamma(R)^{(2)}$ by Example~\ref{gexam1}. Therefore,
$\gamma(\Gamma(R)^{(2)})\geq 2$ by Lemma~\ref{genusbi}.
\\ To finish the proof, let $R\cong R_1\times R_2 \times R_3$.
Assume further that $|R_1|\leq |R_2|\leq |R_3|$ and $|R_3|\geq 3$.
Let $u_1=(1, 0, 0)$, $u_2=(0, 1, 0)$, $u_3=(1, 1, 0)$, $v_1=(0, 0,
1)$, and $v_2=(0, 0, 2)$; then $u_iv_j=0$ for all $i, j$, so that
$K_{2,3}$ is a subgraph of $\Gamma(R)$ it follows that $K_{4,6}$
is a subgraph of $\Gamma(R)^{(2)}$ by Proposition~\ref{gprop}(d).
Therefore, $\gamma(\Gamma(R)^{(2)})\geq 2$ by Lemma~\ref{genusbi},
a contradiction.
\end{proof}

There are some concrete examples for which
$\gamma(\Gamma_I(R))=1$.

\begin{example}
\label{keyexam1} Let $R$ be a finite ring and $I$ be a proper
nonzero ideal of $R$. Then the following hold: \begin{description}
\item{(a)} Suppose that $R/I\cong \mathbb{Z}_2\times
\mathbb{Z}_2\times \mathbb{Z}_2$ and $|I|=2$. Then
$\gamma(\Gamma_I(R))=1$. \item{(b)} Suppose that $R/I\cong
\mathbb{Z}_{16}$ and and $|I|=2$. Then $\gamma(\Gamma_I(R))=1$.
\item{(c)} Suppose that $\Gamma(R/I)=P_3$, $R/I$ is local and
$|I|=3$. Then $\gamma(\Gamma_I(R))=1$. \item{(d)} Suppose that
$R/I\cong \mathbb{Z}_9$ or $R/I\cong
\dfrac{\mathbb{Z}_2[x]}{(x^3)}$ and $|I|=3$. Then
$\gamma(\Gamma_I(R))=1$. \end{description}
\end{example}
\begin{proof}
(a)  By assumption, $\omega(\Gamma(R/I))=3$, therefore
$\Gamma_I(R)$ is not planar by \cite[Theorem~7.2]{r}. By
Remark~\ref{keyrem1}(c), the graph $\Gamma_I(R)$ and the graph
$\Gamma_{0\times 0 \times 0\times \mathbb{Z}_2}(\mathbb{Z}_2\times
\mathbb{Z}_2\times \mathbb{Z}_2\times \mathbb{Z}_2)$ are the same.
Thus, we may assume that $R=\mathbb{Z}_2\times \mathbb{Z}_2\times
\mathbb{Z}_2\times \mathbb{Z}_2$ and $I=0\times 0 \times 0\times
\mathbb{Z}_2$. Let $u_1=(1, 0, 0, 0)$, $u_2=(0, 1, 0, 0), u_3=(0,
0, 1, 0)$, $v_1=(1, 0, 0, 1)$, $v_2=(0, 1, 0, 1), v_3=(0, 0, 1,
1)$, $w_{11}=(0, 1, 1, 0)$, $w_{12}=(0, 1, 1, 1)$, $w_{21}=(1, 0,
1, 0)$, $w_{22}=(1, 0, 1, 1)$, $w_{31}=(1, 1, 0, 0)$ and
$w_{32}=(1, 1, 0, 1)$. Let $G=\Gamma_I(R)$; then $$V(G)=\{u_1,
u_2, u_3, v_1, v_2, v_3, w_{11}, w_{12}, w_{21}, w_{22}, w_{31},
w_{32}\}$$ and
$$\ba{cl} E(G)  = & \{u_iu_j, v_iv_j, u_iv_j, v_iu_j~|~1\leq i<j\leq 3\}\\  &
\cup \{u_iw_{ij}, v_iw_{ij}~|~i=1,2,3,~j=1,2 \}\ea .$$ To draw $G$
on a torus, we first draw a subgraph of $G$ on a torus. Let $G_1$
be the subgraph of $G$ induced by $\{u_1, u_2, u_3, v_1, v_2,
v_3\}$; then Figure 3-1 gives the presentation of $G_1$ on a
torus. To obtain a presentation of $G$, one can simply insert
$\{w_{11}, w_{12}, u_1w_{11}, u_1w_{12}, v_1w_{11}, v_1w_{12}\}$
into the face $F_1$, $\{w_{2j}, u_2w_{2j}, v_2w_{2j}~|~j=1,2\}$
into the face $F_2$ and $\{w_{3j}, u_3w_{3j},
v_3w_{3j}~|~j=1,2\}$ into the face $F_3$.\\

(b) By assumption, $\Gamma(R/I)$ contains a triangle, therefore
$\Gamma_I(R)$ is not planar by \cite[Theorem~7.2]{r}. By
Remark~\ref{keyrem1}(c), the graph $\Gamma_I(R)$ and the graph
$\Gamma_{0\times \mathbb{Z}_2}(\mathbb{Z}_{16}\times
\mathbb{Z}_2)$ are the same. Thus, we may assume that
$R=\mathbb{Z}_{16}\times \mathbb{Z}_2$ and $I=0\times
\mathbb{Z}_2$. Let $u_1=(4, 0)$, $u_2=(8, 0), u_3=(12, 0)$,
$v_1=(4, 1)$, $v_2=(8, 1), v_3=(12, 1)$, $w_1=(2, 0)$, $w_2=(6,
0)$, $w_3=(10, 0)$, $w_4=(14, 0)$, $w_5=(2, 1)$, $w_6=(6, 1)$,
$w_7=(10, 1)$ and $w_8=(14, 1)$. Let $G=\Gamma_I(R)$; then
$$V(G)=\{u_1, u_2, u_3, v_1, v_2, v_3, w_1, \dots, w_8\}$$ and
$$\ba{cl} E(G)  = & \{u_iu_j, v_iv_j, u_iv_j, v_iu_j~|~1\leq i<j\leq 3\}\\ &  \cup\{u_iv_i~|~i=1,2,3\} \\  &
\cup \{u_2w_j, v_2w_j~|~j=1, \dots, 8\}\ea .$$ To draw $G$ on a
torus, we first draw a subgraph of $G$ on a torus. Let $G_1$ be
the subgraph of $G$ induced by $\{u_1, u_2, u_3, v_1, v_2, v_3\}$;
then Figure 3-2 gives the presentation of $G_1$ on a torus. To
obtain a presentation of $G$, one can simply insert $\{w_j,
u_2w_j, v_2w_j~|~j=1, \dots, 8\}$ into the face $F$ of Figure 3-2. \\

\begin{center}
\includegraphics{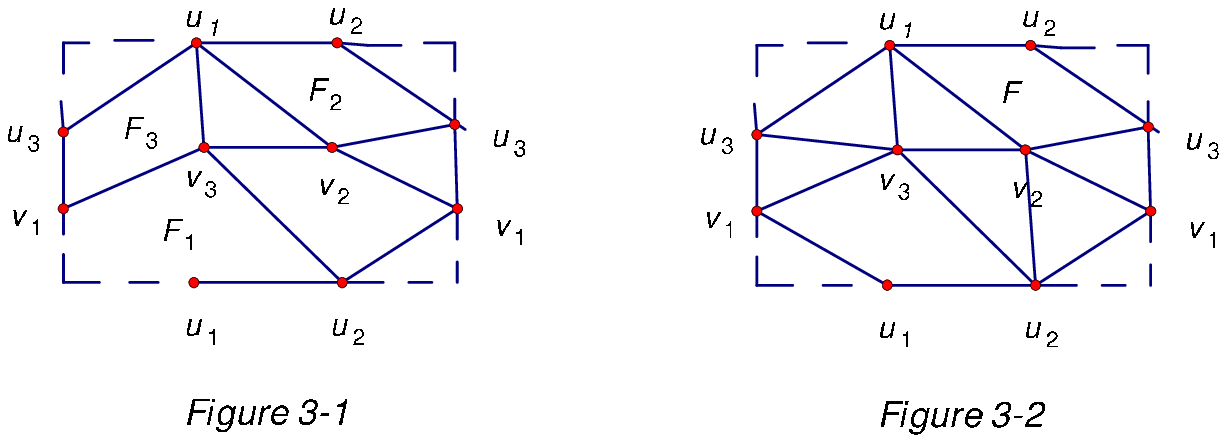}
\end{center}

(c)  By Theorem~\ref{watheo1}, $R/I\cong \mathbb{Z}_8$ or
$R/I\cong \frac{\mathbb{Z}_2[x]}{(x^3)}$. Since the graph
$\Gamma_I(R)$ can be viewed as the graph $\Gamma_{0\times
\mathbb{Z}_3}(R/I\times \mathbb{Z}_3)$ and the graphs
$\Gamma_{0\times \mathbb{Z}_3}(\mathbb{Z}_8\times \mathbb{Z}_3)$
and $\Gamma_{0\times
\mathbb{Z}_3}(\frac{\mathbb{Z}_2[x]}{(x^3)}\times \mathbb{Z}_3)$
are the same, we may assume that $R=\mathbb{Z}_8\times
\mathbb{Z}_3$ and $I=0\times \mathbb{Z}_3$. Let $u_1=(4, 0)$,
$u_2=(4, 1), u_3=(4, 2)$, $v_1=(2, 0)$, $v_2=(2, 1), v_3=(2, 2)$,
$v_4=(6, 0)$, $v_5=(6, 1)$ and $v_6=(6, 2)$. Let $G=\Gamma_I(R)$;
then
$$V(G)=\{u_1, u_2, u_3, v_1, \dots, v_6\}$$ and
$$\ba{cl} E(G)  = & \{u_iv_j~|~1\leq i\leq 3,~1\leq j\leq 6\}\\  &
\cup \{u_iu_j~|~1\leq i<j\leq 3\}\ea .$$  Figure 4 gives the
presentation of $G$ on a torus.
\begin{center}
\includegraphics{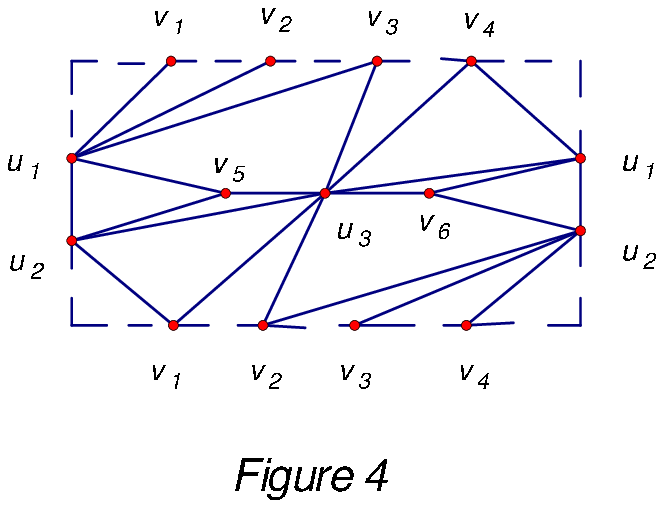}
\end{center}

(d) If $R/I\cong \mathbb{Z}_9$ or $R/I\cong
\dfrac{\mathbb{Z}_2[x]}{(x^3)}$, then $a^2=0$ for every element in
$\Gamma(R/I)$, so that $\Gamma_I(R)=K_6$, it follows that
$\gamma(\Gamma_I(R))=1$.
\end{proof}

To illustrate examples of $\gamma(\Gamma_I(R))\geq 2$, we need a
Lemma.

\begin{lemma} \label{mainlem7} Let $(R, \m)$ be a finite local ring and $I$ be a proper nonzero ideal of $R$.
Suppose that the presentation of $\Gamma(R/I)$ is either $G$ or
$G'$ and $\m/I=\{0, u, v_1, \dots, v_6\}$ such that
$u^2=v_5^2=v_6^2=0$. Then $\gamma(\Gamma_I(R))\geq 2$.

\begin{center}
\includegraphics{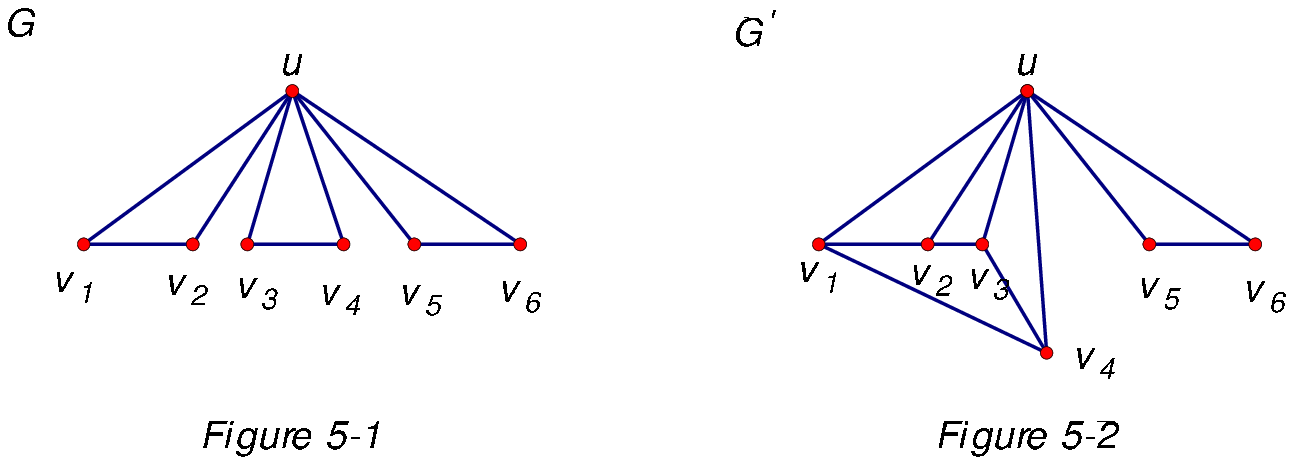}
\end{center}
\end{lemma}

\begin{proof}
We may assume that $|I|=2$ and therefore the graph $\Gamma_I(R$
can be viewed as the graph $\Gamma_{0\times
\mathbb{Z}_2}(R/I\times \mathbb{Z}_2)$ By abuse of notation, let
$u=(u, 0)$ and $v_i=(v_i, 0)$ for $i=1, \dots, 6$. Furthermore,
let $u'=(u, 1)$ and $v_i'=(v_i, 1)$ for $i=1, \dots, 6$. Then
$$V(\Gamma_I(R))=\{u, u', v_1, \dots, v_6, v_1', \dots, v'_6 \}$$
and $$\ba{rl} E(\Gamma_I(R))  \supseteq & \{uu', uv_i, uv'_i,
u'v_i, u'v'_i~|~i=1, \dots, 6\}\\ &  \cup \{v_5v_5', v_5v_6,
v_5v_6', v_5'v_6, v_5'v_6', v_6v_6'\}\\ &  \cup \{v_iv_{i+1},
v_iv_{i+1}', v_i'v_{i+1}, v_i'v_{i+1}'~|~i=1, 3\}  \ea .
$$
Let $H$ be the subgraph of $\Gamma_I(R)$ such that
$V(H)=V(\Gamma_I(R))$ and $$\ba{rl} E(H)= & \{uu', uv_i, uv'_i,
u'v_i, u'v'_i~|~i=1, \dots, 6\}\\ &  \cup \{v_5v_5', v_5v_6,
v_5v_6', v_5'v_6, v_5'v_6', v_6v_6'\}\\ &  \cup \{v_iv_{i+1},
v_iv_{i+1}', v_i'v_{i+1}, v_i'v_{i+1}'~|~i=1, 3\} \ea .
$$ Let $H_1=H-\{v_5, v_5', v_6, v_6'\}$
and $H_2$ be the subgraph induced by $\{v_5, v_5', v_6, v_6'\}$;
then there are edges $\{uv_5, uv_5', uv_6, uv_6'\}$ of $H$ join
$H_2$ to $H_1$ and $H_2=K_4$, it follows that
$\gamma(H)>\gamma(H_1)$ by Lemma~\ref{keylem}. Therefore, to
finish the proof it is enough to show that $H_1$ is not planar.
For this, observe that $H_1$ has $10$ vertices and $25$ edges. If
$H_1$ can be drawn in a plane, then $H_1$ has $17$ faces by Euler
characteristic formula. However, the edges and faced satisfies
$50\geq 17\cdot 3$, a contradiction. Thus, $\gamma(H_1)\geq 1$.
\end{proof}

\begin{example}
\label{keyexam2} Let $R$ be a finite ring and $I$ be a proper
nonzero ideal of $R$.\\ (i) If $R/I$ is one of the following local
rings: (a)~${\mathbb{Z}_4[x, y]}/{(x^2, y^2, xy-2)}$,
(b)~${\mathbb{Z}_2[x, y]}/{(x^2, y^2)}$,
(c)~${\mathbb{Z}_4[x]}/{(x^2)}$, then $\gamma(\Gamma_I(R))\geq 2$. \\
(ii) If $R/I$ is one of the following local rings:
(a)~${\mathbb{Z}_2[x, y]}/{(x^3, xy, y^2-x^2)}$,
(b)~${\mathbb{Z}_4[x]}/{(x^3, x^2-2x)}$,
(c)~${\mathbb{Z}_4[x,y]}/{(x^3, x^2-2, xy, y^2-2)}$,
(d)~${\mathbb{Z}_8[x]}/{(x^2-4, 2x)}$, then
$\gamma(\Gamma_I(R))\geq 2$.
\end{example}
\begin{proof}
We may assume that $|I|=2$. \\ (i) If $R/I$ is one of the rings in
(ii), then the graph of $\Gamma(R/I)$ is $G'$ in Figure 5-1 by
\cite[Example~2.7]{cw}. Hence by Lemma~\ref{mainlem7}, we need to
find $3$ elements in each ring with nilpotency 2. If
$R={\mathbb{Z}_4[x, y]}/{(x^2, y^2, xy-2)}$, then $u=\bar 2$,
$v_5=\bar x+\bar y$ and $v_6=\bar x+\bar y+\bar 2$ are the
required 3 elements. If $R={\mathbb{Z}_2[x, y]}/{(x^2, y^2)}$,
then $u=\bar x\bar y$, $v_5=\bar x+\bar y$ and $v_6=\bar x+\bar
y+\bar x\bar y$ are the required 3 elements. If
$R={\mathbb{Z}_4[x]}/{(x^2)}$, then $u=\bar{2x}$, $v_5=\bar 2+\bar
x$ and $v_6=\bar 2+\bar{3x}$ are
the required 3 elements.\\
(ii) If $R/I$ is one of the ring in (i), then the graph of
$\Gamma(R/I)$ is $G'$ in Figure 5-2 by \cite[Example~2.6]{cw}.
Hence by Lemma~\ref{mainlem7}, we need to find $3$ elements in
each ring with nilpotency 2. If $R={\mathbb{Z}_2[x, y]}/{(x^3,xy,
y^2-x^2)}$, then $u=\bar x^2$, $v_5=\bar x+\bar y$, and $v_6=\bar
x+\bar y+\bar x^2$ are the required 3 elements. If
$R={\mathbb{Z}_4[x]}/{(x^3, x^2-2x)}$, then $u= \bar{2x}$,
$v_5=\bar 2$ and $v_6=\bar 2+ \bar{2x}$ are the required 3
elements. If $R={\mathbb{Z}_4[x,y]}/{(x^3, x^2-2, xy, y^2-2)}$,
then $u=\bar 2$, $v_5=\bar x+\bar y$ and $v_6=\bar x+\bar y+\bar
2$ are the required 3 elements. If $R={\mathbb{Z}_8[x]}/{(x^2-4,
2x)}$, then $u=\bar 4$, $v_5=\bar 2+\bar x$ and $v_6=\bar 6+\bar
x$ are the required 3 elements.
\end{proof}

There are some sufficient conditions for which
$\gamma(\Gamma_I(R))\geq 2$.
\begin{lemma}
\label{mainlem8} Let $R$ be a finite ring and $I$ be a proper
nonzero ideal of $R$. Suppose that $\gamma(\Gamma(R/I))\geq 2$.
Then $\gamma(\Gamma_I(R))\geq 2$.
\end{lemma}
\begin{proof}
The conclusion follows  as $$2\leq \gamma(\Gamma(R/I))\leq
\gamma(\Gamma(R/I)^{(2)})\leq \gamma(\Gamma_I(R)).$$
\end{proof}

\begin{lemma}
\label{mainlem9} Let $R$ be a finite ring and $I$ be a proper
nonzero ideal of $R$. Suppose the following hold:
\begin{description}
\item{(a)} $R/I$ is local with unique maximal ideal $\m/I$.
\item{(b)} $|R/\m|=2$. \item{(c)} $\gamma(\Gamma(R/I))=1$.
\end{description}
Then $\gamma(\Gamma_I(R))\geq 2$.
\end{lemma}
\begin{proof}
By (a) and (c), $R/I$ is isomorphic to one of the rings in
Theorem~\ref{watheo2}. However, if we take (b) into account, then
$K_{3,4}$ is a subgraph of every ring in Theorem~\ref{watheo2}
with (b) holds. Thus, $2\leq \gamma(K_{6,8})\leq
\gamma(\Gamma(R/I)^{(2)})\leq \gamma(\Gamma_I(R))$.
\end{proof}

\begin{theo}\label{keytheo1}
Let $R$ be a finite ring and $I$ be a proper nonzero ideal of $R$.
Suppose that $\omega(\Gamma(R/I))\geq 4$ or
$\gamma(\Gamma(R/I))\geq 1$. Then $\gamma(\Gamma_I(R))\geq 2$.
\end{theo}
\begin{proof}
Since $R$ is finite, $R/I\cong R_1\times \cdots \times R_k$, where
$R_i$ is a finite local ring for every $i$. If $k\geq 4$, then by
Lemma~\ref{mainlem6} and the facts that $|I|\geq 2$ and
$\Gamma((R/I)^{(2)})$ is a subgraph of $\Gamma_I(R)$, we see that
$\gamma(\Gamma_I(R))\geq 2$. Therefore, we may assume that $k\leq
2$ or $R/I\cong \mathbb{Z}_2\times \mathbb{Z}_2\times
\mathbb{Z}_2$. However, $R\cong \mathbb{Z}_2\times
\mathbb{Z}_2\times \mathbb{Z}_2$ implies that
$\omega(\Gamma(R/I))=3$ and $\Gamma(R/I)$ is planar, a
contradiction. Thus, $k\leq 2$.\par Now, suppose that $R/I$ is not
local and $R/I\cong R_1\times R_2$ with $|R_1|\leq |R_2|$.  If
$|R_1|\geq 3$ and $|R_2|\geq 4$, then $K_{2,3}$ is a subgraph of
$\Gamma(R/I)$, so that $\gamma(\Gamma(R/I)^{(2)})\geq 2$ by
Proposition~\ref{gprop}(d), it follows that
$\gamma(\Gamma_I(R))\geq 2$. If $|R_1|=|R_2|=3$, then $R/I\cong
\mathbb{Z}_3\times \mathbb{Z}_3$, which contradicts to the
assumptions that $\omega(\Gamma(R/I))\geq 4$ or
$\gamma(\Gamma(R/I))\geq 1$. So, we may assume that $R_1\cong
\mathbb{Z}_2$. Since $\gamma(\Gamma(\mathbb{Z}_2\times R_2))\geq
1$ or $\omega(\Gamma(\mathbb{Z}_2\times R_2))\geq 4$, by
Lemma~\ref{mainlem4}, $R_2$ is neither a field nor satisfies
$|\Gamma(R_2)|=1$. Therefore $|\Gamma(R_2)|\geq 2$. However, by
Lemma~\ref{mainlem4} again $K_{2,3}$ is a subgraph of
$\Gamma(R/I)$, it follows that $\gamma(\Gamma(R/I)^{(2)})\geq 2$
by Proposition~\ref{gprop}(d). Hence, $\gamma(\Gamma_I(R))\geq
2$.\par Finally we assume that $R/I$ is local with unique maximal
ideal $\m/I$.\\  Case~1.  $(\m/I)^2=0$. In this case,
$\omega(\Gamma(R/I))\geq 4$ or $\gamma(\Gamma(R/I))\geq 1$ implies
that there are  $4$ distinct nonzero zero-divisors $u_1, \dots,
u_4$ in $R/I$ such that $u_iu_j=0$ for all $i, j$, so that $K_8$
is a subgraph of $\Gamma(R/I)^{(2)}$, it follows that
$\gamma(\Gamma_I(R))\geq 2$. \\ Case~2.  $(\m/I)^2\neq 0$. In this
case, observe $\Gamma(R/I)^{(2)}$ is a subgraph of $\Gamma_I(R)$.
Hence, by Lemma~\ref{mainlem2}, we may assume that $|R/\m|=2$ as
$|R/\m|\geq 3$ implies that $\gamma(\Gamma_I(R))\geq 2$.
Therefore, $|\Gamma(R/I)|=2^n-1$ for some positive integer $n$.
Observe that $|\Gamma(R/I)|$ is 1 or 3 are not possible. Moreover,
the assumption that $\gamma(\Gamma(R/I))\geq 1$ implies that
$\gamma(\Gamma_I(R))\geq 2$ by Lemma~\ref{mainlem8} and
Lemma~\ref{mainlem9}. Therefore, we may assume finally that
$\omega(\Gamma(R/I))\geq 4$ and $\Gamma(R/I)$ is planar. However,
if so, then by Theorem~\ref{watheo1}, $|\Gamma(R/I)|=7$ and $R/I$
is isomorphic to $\mathbb{Z}_{16}$ or isomorphic to one of the
rings in Example~\ref{keyexam2}, so that $\omega(\Gamma(S))=3$, a
contradiction. The proof is now complete.
\end{proof}

\begin{theo}
\label{keytheo2} Let $R$ be a finite ring and $I$ be a proper
nonzero ideal of $R$ which is not a prime ideal. Suppose that
$\omega(\Gamma(R/I))\leq 2$. Then $\gamma(\Gamma_I(R))\leq 1$ if
and only if one of the following holds:
\begin{description}
\item{(a)}  $R/I\cong \mathbb{Z}_3\times \mathbb{Z}_3$ and
$|I|\leq 2$. \item{(b)} $R/I\cong \mathbb{Z}_2\times \mathbb{Z}_2$
and $|I|\leq 4$. \item{(c)}  $R/I\cong \mathbb{Z}_2\times
\mathbb{Z}_3$ and $|I|\leq 3$.\item{(d)} $R/I\cong
\mathbb{Z}_2\times \mathbb{Z}_4$ and $|I|\leq 2$. \item{(e)}
$R/I\cong \mathbb{Z}_2\times \frac{\mathbb{Z}_2[x]}{(x^2)}$ and
$|I|\leq 2$. \item{(f)} $R/I\cong \mathbb{Z}_2\times \mathbb{F}_q$
with $q\geq 4$ and $|I|\leq 2$.  \item{(g)} $R/I=\mathbb{Z}_4$ and
$|I|\leq 7$. \item{(h)} $R/I=\frac{\mathbb{Z}_2[x]}{(x^2)}$ and
$|I|\leq 7$. \item{(i)} $R/I=\mathbb{Z}_9$ and $|I|\leq 3$.
\item{(j)} $R/I=\frac{\mathbb{Z}_3[x]}{(x^2)}$ and $|I|\leq 3$.
\item{(k)} $R/I=\mathbb{Z}_8$ and $|I|\leq 3$.\item{(l)}
$R/I=\frac{\mathbb{Z}_2[x]}{(x^3)}$ and $|I|\leq 3$.\item{(m)}
$R/I=\frac{\mathbb{Z}_4[x]}{(x^2-2, x^3)}$ and $|I|\leq 3$.
\end{description}
\end{theo}

\begin{proof}
($\Rightarrow$) Assume that $\gamma(\Gamma_I(R))\leq 1$. From the
first paragraph of the proof of Theorem~\ref{keytheo1}, we may
assume that $R/I=\mathbb{Z}_2\times \mathbb{Z}_2\times
\mathbb{Z}_2$ or $R/I$ is local or $R/I$ is a product of two
finite local rings. However, $R/I=\mathbb{Z}_2\times
\mathbb{Z}_2\times \mathbb{Z}_2$ implies that
$\omega(\Gamma(R/I))=3$, a contradiction. Thus, $R/I$ is local or
a product of two finite local rings.  \par Now, suppose that $R/I$
is not local and $R/I\cong R_1\times R_2$ with $|R_1|\leq |R_2|$.
From the second paragraph of the proof of Theorem~\ref{keytheo1},
we may assume that $R/I\cong \mathbb{Z}_3\times \mathbb{Z}_3$ or
$R_1\cong \mathbb{Z}_2$. If $R/I\cong \mathbb{Z}_3\times
\mathbb{Z}_3$, then $\Gamma(R/I)=C_4=K_{2,2}$ and $R/I$ has no
nilpotent, so that $|I|\leq 2$. Now, we assume that
$R_1=\mathbb{Z}_2$. By Lemma~\ref{mainlem4}, $R_2$ is either a
field or satisfies $|\Gamma(R_2)|=1$. It follows that $R/I$ is
isomorphic to one of the following rings: $\mathbb{Z}_2\times
\mathbb{Z}_2$, $\mathbb{Z}_2\times \mathbb{Z}_3$,
$\mathbb{Z}_2\times \mathbb{Z}_4$, $\mathbb{Z}_2\times
\frac{\mathbb{Z}_2[x]}{(x^2)}$ and $\mathbb{Z}_2\times
\mathbb{F}_q$, where $q\geq 4$. \\ If $R/I\cong \mathbb{Z}_2\times
\mathbb{Z}_2$, then $\Gamma(R/I)=K_2$. Since $R/I$ has no
nilpotent, $\Gamma_I(R)=K_{t,t}$ if $|I|=t$, it follows that
$|I|\leq 4$ by Lemma~\ref{genusbi}.\\ If $R/I\cong
\mathbb{Z}_2\times \mathbb{Z}_3$, then $\Gamma(R/I)=K_{1,2}$.
Since $R/I$ has no nilpotent, $\Gamma_I(R)=K_{t,2t}$ if $|I|=t$,
it follows that $|I|\leq 3$ by Lemma~\ref{genusbi}.\\ If $R/I$ is
isomorphic to $\mathbb{Z}_2\times \mathbb{Z}_4$ or
$\mathbb{Z}_2\times \frac{\mathbb{Z}_2[x]}{(x^2)}$ or
$\mathbb{Z}_2\times \mathbb{F}_q$ with $q\geq 4$, then $K_{1,3}$
is a subgraph of $\Gamma(R/I)$, it follows that $|I|\leq 2$ by
Lemma~\ref{genusbi}.
\par From the above, we can finally assume that $R/I$ is local with unique
maximal ideal $\m/I$. By Theorem~\ref{keytheo1}, $\Gamma(R/I)$ is
planar. So we only consider the rings in Theorem~\ref{watheo1}.\\
Case~1. $(\m/I)^2=0$. In this case, $\Gamma(R/I)\cong K_t$ if
$|(\m/I)|=t+1$, so that $\Gamma(R/I)$ is $K_1$ or $K_2$ as
$\omega(\Gamma(R/I))\leq 2$, it follows that $R/I$ is isomorphic
to the following rings: $\mathbb{Z}_4$,
$\frac{\mathbb{Z}_2[x]}{(x^2)}$, $\mathbb{Z}_9$ and
$\frac{\mathbb{Z}_3[x]}{(x^2)}$. \\ If $R/I$ is isomorphic to
$\mathbb{Z}_4$ or $\frac{\mathbb{Z}_2[x]}{(x^2)}$, then
$\Gamma_I(R)=K_t$ if $|I|=t$, it follows that $|I|\leq 7$ by
Lemma~\ref{genuskn}.\\ If $R/I$ is isomorphic to $\mathbb{Z}_9$ or
$\frac{\mathbb{Z}_3[x]}{(x^2)}$, then $\Gamma_I(R)=K_{2t}$ if
$|I|=t$, it follows that $|I|\leq 3$ by
Lemma~\ref{genuskn}.\\
Case~2. $(\m/I)^2\neq 0$. In this case,  $|R/\m|=2$ as $|R/\m|\geq
3$ implies that $\gamma(\Gamma_I(R))\geq 2$ ($\Gamma(R/I)^{(2)}$
is a subgraph of $\Gamma_I(R)$) by Lemma~\ref{mainlem2}. Moreover,
by Theorem~\ref{keytheo1}, $\Gamma(R/I)$ is planar. Hence by
Theorem~\ref{watheo1}, $|\Gamma(R/I)|=1, 3, 7$. From the above, we
see that $|\Gamma(R/I)|=3, 7$. However, if $|\Gamma(R/I)|=7$, then
$R/I$ is isomorphic to $\mathbb{Z}_{16}$ or isomorphic to one of
the rings in Example~\ref{keyexam2}, so that $\Gamma(R/I)$
contains a triangle, a contradiction. Thus, $|\Gamma(R/I)|=3$ and
$|\Gamma(R/I)$ is $P_3$ as $\Gamma(R/I)$ contains no triangles and
is connected. Therefore, $R/I$ is isomorphic to the following
rings: $\mathbb{Z}_8$, $\frac{\mathbb{Z}_2[x]}{(x^3)}$ and
$\frac{\mathbb{Z}_4[x]}{(x^2-2, x^3)}$. In either cases,
$K_{t,2t}$ is a subgraph of $\Gamma_I(R)$ if $|I|=t$, therefore
$t\leq 3$.

($\Leftarrow$) Assume that $R/I\cong \mathbb{Z}_3\times
\mathbb{Z}_3$ and $|I|\leq 2$. Then $\Gamma(R/I)=K_{2,2}$ and
$R/I$ has no nilpotent. It follows that $\Gamma_I(R)=K_{2t,2t}$
with $|I|=t\leq 2$. Therefore, $\gamma(\Gamma(R/I))\leq 1$ by Lemma~
\ref{genusbi}.\\
Assume that $R/I\cong \mathbb{Z}_2\times \mathbb{Z}_2$ and
$|I|\leq 4$. Then $\Gamma(R/I)=K_2$ and $R/I$ has no nilpotent. It
follows that $\Gamma_I(R)=K_{2t}$ with $|I|=t\leq 4$. Therefore,
$\gamma(\Gamma(R/I))\leq 1$ by Lemma~ \ref{genusbi}.\\ Assume that
$R/I\cong \mathbb{Z}_2\times \mathbb{Z}_3$ and $|I|\leq 3$. Then
$\Gamma(R/I)=K_{1,2}$ and $R/I$ has no nilpotent. It follows that
$\Gamma_I(R)=K_{t,2t}$ with $|I|=t\leq 3$. Therefore,
$\gamma(\Gamma(R/I))\leq 1$ by Lemma~ \ref{genusbi}.\\ Assume that
$R\cong \mathbb{Z}_2\times \mathbb{Z}_4$ or $R\cong
\mathbb{Z}_2\times \frac{\mathbb{Z}_2[x]}{(x^2)}$ or $R/I\cong
\mathbb{Z}_2\times \mathbb{F}_q$ with $q\geq 4$ and $|I|\leq 2$.
Then $\Gamma_I(R)$ is planar by \cite[Theorem~7.2]{r}.\\
Assume that $R\cong \mathbb{Z}_4$ or $R\cong
\frac{\mathbb{Z}_2[x]}{(x^2)}$. In either cases, $\Gamma_I(R)=K_t$
if $|I|=t$. Therefore if $|I|\leq 7$, then
$\gamma(\Gamma_I(R))\leq 1$ by Lemma~\ref{genuskn}. \\
Assume that $R/I\cong \mathbb{Z}_9$ or $R/I\cong
\mathbb{Z}_3[x]/(x^2)$. Then for every zero-divisor $a\in R/I$
satisfies $a^2=0$. Therefore, $\Gamma_I(R)=K_{2t}$. if $|I|=t$.
Thus, if $|I|\leq 3$, then $\gamma(\Gamma_I(R))\leq 1$ by
Lemma~\ref{genuskn}. \\
Assume that $R/I$ is isomorphic to the following rings:
$\mathbb{Z}_8$, $\frac{\mathbb{Z}_2[x]}{(x^3)}$ and
$\frac{\mathbb{Z}_4[x]}{(x^2-2, x^3)}$. In either cases, if
$|I|\leq 3$ then $\gamma(\Gamma_I(R))\leq 1$ by
Example~\ref{keyexam1}(c).
\end{proof}

\begin{theo}
\label{keytheo3} Let $R$ be a finite ring and $I$ be a proper
nonzero ideal of $R$. Suppose that $\omega(\Gamma(R/I))=3$. Then
$\gamma(\Gamma_I(R))=1$ if and only if $|I|=2$ and $R/I$ is
isomorphic to one of the following:
\begin{description}
\item{(a)} $\mathbb{Z}_2\times \mathbb{Z}_2\times \mathbb{Z}_2$.
\item{(b)} $\mathbb{Z}_{16}$. \item{(c)}
$\dfrac{\mathbb{Z}_2[x,y]}{(x^2, xy, y^2)}$. \item{(d)}
$\dfrac{\mathbb{Z}_4[x]}{(2x, x^2)}$. \item{(e)}
$\dfrac{\mathbb{F}_4[x]}{(x^2)}$.\item{(f)}
$\dfrac{\mathbb{Z}_4[x]}{(x^2+x+1)}$.
\end{description}
\end{theo}
\begin{proof}
($\Rightarrow$) Assume that $\gamma(\Gamma_I(R))=1$. From the
first paragraph of the proof of Theorem~\ref{keytheo1}, we may
assume
that $R/I\cong \mathbb{Z}_2\times \mathbb{Z}_2\times \mathbb{Z}_2$ or $R/I$ is local or $R/I$ is a product of two finite local rings. \\
Assume that $R/I\cong \mathbb{Z}_2\times \mathbb{Z}_2\times
\mathbb{Z}_2$. If $|I|\geq 3$, then $\gamma(\Gamma(R/I)^{(3)})\geq
2$ as $K_{1,3}$ is a subgraph of $\Gamma(R/I)$, it follows that
$\gamma(\Gamma_I(R))\geq 2$ . Moreover, If $|I|=2$, then
$\gamma(\Gamma_I(R))=1$ by Example~\ref{keyexam1}(a).\\
Assume that $R/I$ is not local and $R/I\cong R_1\times R_2$ with
$|R_1|\leq |R_2|$. From the second paragraph of the proof of
Theorem~\ref{keytheo1}, we may assume that $R/I\cong
\mathbb{Z}_3\times \mathbb{Z}_3$ or $R_1\cong \mathbb{Z}_2$.
However, if $R/I\cong \mathbb{Z}_3\times \mathbb{Z}_3$, then
$\omega(\Gamma(R/I))=2$, a contradiction. Thus, $R_1\cong
\mathbb{Z}_2$. Moreover, if $R_2$ is a field or $\Gamma(R_2)$ is a
point then $\Gamma(R/I)$ contains no cycles by
Lemma~\ref{mainlem4} which contradicts to the assumption that
$\omega(\Gamma(R/I))=3$. If $|\Gamma(R_2)|\geq 2$, then $K_{2,3}$
is a subgraph of $\Gamma(R/I)$, it follows that
$\gamma(\Gamma_I(R))\geq 2$, a contradiction again.
\par From the above, we can finally assume that $R/I$ is local with unique
maximal ideal $\m/I$. Observe that the assumption $\omega(\Gamma(R/I))=3$ implies that $|(\m/I)|\geq 4$.\\
Case~1. $(\m/I)^2=0$. In this case, $\Gamma(R/I)\cong K_t$ if
$|(\m/I)|=t+1$, so that $\Gamma(R/I)\cong K_3$ as
$\omega(\Gamma(R/I))=3$, it follows that $R/I$ is isomorphic to
$\dfrac{\mathbb{Z}_2[x,y]}{(x^2, xy, y^2)}$ or
$\dfrac{\mathbb{Z}_4[x]}{(2x, x^2)}$ or
$\dfrac{\mathbb{F}_4[x]}{(x^2)}$ or
$\dfrac{\mathbb{Z}_4[x]}{(x^2+x+1)}$ by Lemma~\ref{mainlem5}.
Moreover, every element $a\in R/I$ satisfies $a^2=0$ Thus,
$\Gamma_I(R)=K_{3t}$ if
$|I|=t$. It follows that $\gamma(\Gamma_I(R))=1$ if and only if $|I|=2$. \\
Case~2. $(\m/I)^2\neq 0$. In this case,  $|R/\m|=2$ as $|R/\m|\geq
3$ implies that $\gamma(\Gamma_I(R))\geq 2$ ($\Gamma(R/I)^{(2)}$
is a subgraph of $\Gamma_I(R)$) by Lemma~\ref{mainlem2}. Moreover,
by Theorem~\ref{keytheo1}, $\Gamma(R/I)$ is planar. Hence by
Theorem~\ref{watheo1}, $|\Gamma(R/I)|=1, 3, 7$. From the above, we
may assume that $|\Gamma(R/I)|=7$. However, if so $R/I$ is
isomorphic to $\mathbb{Z}_{16}$ or isomorphic to one of the rings
in Example~\ref{keyexam2}. Since $\gamma(\Gamma_I(R))\geq 2$ if
$R/I$ is isomorphic to one of the rings in Example~\ref{keyexam2},
we see that $R/I\cong \mathbb{Z}_{16}$. If $|I|\geq 3$, then
$\gamma(\Gamma_I(R))\geq 2$ as $K_{1,3}$ is a subgraph of
$\Gamma(R/I)$. Thus $|I|=2$.\\
($\Leftarrow$) Assume that $R/I$ is isomorphic to one of the
following $2$ rings: $\mathbb{Z}_2\times \mathbb{Z}_2\times
\mathbb{Z}_2$ and $\mathbb{Z}_{16}$. If
$|I|=2$, then $\gamma(\Gamma_I(R))=1$ by Example~\ref{keyexam1}.\\
Assume that $R/I$ is isomorphic to one of the following $4$ rings:
$\dfrac{\mathbb{Z}_2[x,y]}{(x^2, xy, y^2)}$,
$\dfrac{\mathbb{Z}_4[x]}{(2x, x^2)}$,
$\dfrac{\mathbb{F}_4[x]}{(x^2)}$ and
$\dfrac{\mathbb{Z}_4[x]}{(x^2+x+1)}$ . If $|I|=2$, then
$\gamma(\Gamma_I(R))=1$ as $\Gamma_I(R)=K_6$.
\end{proof}


\begin{thebibliography}{1}


\bibitem{amy}
S. Akbari, H.R. Maimani and S. Yassemi,
\newblock {\em When a zero-divisor graph is planar or a complete $r$-partite graph},
\newblock J. Algebra {\bf 270} (2003), 169--180.


\bibitem{am}
S. Akbari and A. Mohammadian,
\newblock {\em On the zero-divisor graph of a commutative ring},
\newblock J. Algebra {\bf 274} (2004), 847--855.




\bibitem{afll}
D.F. Anderson, A. Frazier, A. Lauve and P.S. Livingston,
\newblock {\em The zero-divisor graph of a commutative ring, II},
\newblock Lecture Notes in Pure and Appl. Math. {\bf 220} (2001), 61--72.

\bibitem{al}
D.F. Anderson and P.S. Livingston,
\newblock {\em The zero-divisor graph of a commutative ring},
\newblock J. Algebra {\bf 217} (1999), 434--447.


\bibitem{cw}
H.-J Chiang-Hsieh and H.-J, Wang,
\newblock {\em Commutative rings with toroidal zero-divisor graphs"},
\newblock preprint.



\bibitem{mpy}
H.R. Maimani, M.R. Pournaki and S. Yassemi,
\newblock {\em Zero-divisor graph with respect to an ideal},
\newblock Comm. Algebra {\bf 34} (2006), no. 3, 923--929.






\bibitem{r}
S. Redmond,
\newblock {\em An ideal-based zero-divisor graph of a commutative ring},
\newblock  Comm. Algebra {\bf 31} (2003), no. 9, 4425--4443.





\bibitem{ry}
G. Ringel and J.W.T.Youngs,
\newblock {\em Solution of the Heawood map-coloring problem},
\newblock Proc. Nat. Acad. Sci. USA {\bf 60} (1968), 438--445.

\bibitem{s}
N.O. Smith,
\newblock {\em Planar zero-divisor graphs},
\newblock International J. Commutative Rings {\bf 2} (2003), no. 4, 177-188.

\bibitem{wa}
H.-J, Wang,
\newblock {\em Zero-divisor graphs of genus one},
\newblock J. Algebra {\bf 304} (2006), no. 2, 666--678.




\end{thebibliography}
\end{document}